%--------- PUBLE.TEX ----------

%------Commenced on 11 November, 2005

%------- draft: December 5, 2005

\def\caption{${\underline{\hbox{\sbf  Haifa Public Talk }}\atop
\hbox{\sbf December 2005}}$}

\magnification=\magstep1

%%%%%%%%%%%%%%% PLAIN TeX %%%%%%%%%%%%%%%

%-----Paper Format---------
\hsize=14.2cm
\vsize=20cm
\hoffset=-0.4cm
\voffset=0.5cm

%-----AMS Symbols-------
\input amssym.def
\input amssym.tex

%-----Special Fonts-----
\font\srm=cmr8
\font\sbf=cmbx8
\font\csc=cmcsc10
\font\scsc=cmcsc10 at 8pt
\font\bcsc=cmcsc10 at 12pt
\font\title=cmr12 at 15pt

\font\teneusm=eusm10
\font\seveneusm=eusm7
\font\fiveeusm=eusm5
\newfam\eusmfam

\textfont\eusmfam=\teneusm
\scriptfont\eusmfam=\seveneusm
\scriptscriptfont\eusmfam=\fiveeusm

\font\teneufm=eufm10
\font\seveneufm=eufm7
\font\fiveeufm=eufm5
\newfam\eufmfam

\textfont\eufmfam=\teneufm
\scriptfont\eufmfam=\seveneufm
\scriptscriptfont\eufmfam=\fiveeufm

%----Small Macros----
\def\varGamma{{\mit \Gamma}}
\def\Re{{\rm Re}\,}

\def\txt#1{{\textstyle{#1}}}

\def\r#1{{\rm #1}}
\def\B#1{{\Bbb #1}}

%----Headers-------
\nopagenumbers
\def\rightheadline{\hfil{\srm Prime Numbers -- Your Gems}
\hfil\tenrm\folio}
\def\leftheadline{\tenrm\folio\hfil{\scsc Y. Motohashi}\hfil}
\def\emptyheadline{\hfil}
\headline{\ifnum\pageno=1 \emptyheadline\else
\ifodd\pageno \rightheadline \else \leftheadline\fi\fi}
%-------Text---------
\vglue 0.7cm
\centerline{\title Prime Numbers -- Your Gems}
\vskip 1cm
\centerline{\bcsc By Yoichi Motohashi}
\vskip 1cm
\noindent
\item{} {\csc Contents}
\item{1.} Introduction
\item{2.} Natural Numbers
\item{3.} Prime Numbers
\item{4.} Sieve Methods
\item{5.} Riemann Hypothesis
\vskip 1cm
\noindent
{\csc Prime} numbers or primes are man's eternal treasures that have been
cherished for several millennia, until today. As their academic ancestors in
ancient Mesopotamia, many mathematicians are
still trying hard to see primes better. It is
my great luck  that I have been able to lead a life of such a mathematician for some
forty years without a single moment of being wearied. I shall relate here a part of my
impressions that I have gathered in a corner of my mind through my own research and
excursions outside my profession.  Mathematicians savour poignant special and surmise
grand general, whence I  shall likewise indulge in extrapolations via  my 
readings, perhaps to a somewhat exceeding extent. I shall be more than rewarded if any
of you shares my reflections and sentiments. 
 \smallskip
This is in essence a translation of my Japanese article that was prepared for
my public talk at the general assembly of the Mathematical Society
of Japan, Spring 2005, and is now available on the home page of MSJ. At this
opportunity I have made substantial revisions, and changed the title into a more
appropriate one.  The use of mathematical formulas and deductions 
is restricted to a minimum, excepting the final section. One may start with
the last three sections, which are more mathematical, 
and turn to the first two, which are more bookish. 
My choice from grammatical genders is made for convenience, save for obvious
situations.
\bigskip
\noindent
{\bf 1. Introduction.} The title of my MSJ public talk was `On the wing of
primes', which bears an obvious borrowing from a piece of F. Mendelssohn;
`songs' was replaced by `primes'. Songs are often equivalent to poems,
specifically in the Japanese tradition. There are a variety of anthologies from
ancient Japan, many of which were compiled by imperial edicts with the earthy aim to
charge the power with grandeur, but only poets therein attained to the ever 
nobility. One of them is  the Kokin-waka-shu (905 AD), and what made it a
man's treasure are not only the kaleidoscopic collection of  those 1100
poems `waka' but also its preface, apparently due to one of the editors, a
highly celebrated poet. It has an
impressive passage: ``That may move the universe without any exertion of power is
poem''. This must be the very element of poems or songs. Indeed, without any exertion
of mechanical power, a tiny song might alter our view of the world and
the firmament. A poet renders salient and subtle  phenomena that are supported by 
invisible existences filling our universe. The universe is perhaps the dual 
or the reflection of the
whole of poems ever created. Later we shall see in a distance that a similar duality
exists in the very fundamentals of the researches on primes.
\smallskip
Primes are rigid existences at the bottom
of the human intelligence as gems and petals of cherry blossoms. Mathematicians
have ever been tempted to lift the veil enveloping primes. 
Sometimes mathematics is
regarded as a precision science that has contributed to the
realization of our modern lifeway. However, mathematicians themselves do not see their
researches as anything that might bring about a further efficient
civilization. Mathematicians ever drift, courting decent theorems. A theorem is
a song.  Japanese mathematicians of the Edo period
(c.\ 1600--1867 AD) enshrined their discoveries and difficult proofs to Shinto
shrines, with gorgeous frameworks; the gods might have been a little bit perplexed
with mathematical formulas and diagrams. 
\smallskip
On the other hand `wing' in the former title suggests that we have been
deepening our thoughts relying on works due to former generations. This is
another salient point to be espied in my discussion. 
Thus, history or the power of letters is also a theme. One may see
again a similitude between mathematics and the whole of poems. One who
brought forth a new idea is ever revered {\it per se\/}. 
The Man-yo-shu (4516 poems, c.\ 430--760 AD), which Japanese often liken
to holy scrolls, enlightens mathematicians as well. The editor a great poet himself
did  painstaking searches for the names of poets and the moments and venues 
of composition.
Thus, even the poems of some rebels against the purple, who were
eventually executed, as well as even bleeding laments of soldiers who were separated
from their families in the remote northern countryside and sent to the front on the
shores of the west islands against the invaders from the continent are mostly
annotated if succinct. Near the end of the preface of the Kokin-waka-shu, the editor
writes, ``The soul of a poet remains if only these letters do.'' This reminds me
of the maxim `Verba volant. Scripta manent.'  
\smallskip
Moving to a different topic, I observe that both the presentation of one's work at
a conference and the teaching at a class room belong to the theatrical
art. Namely, they are reconstructions of something created already. I do not mean
that a theatrical performance is not a creative art; on the contrary, it can indeed
be a creation of something quite new, and audiences come to witness such moments. I
wish you to notice that an important characteristic of a mathematician is his ability
of reconstructing his own or others' discoveries, which is somewhat similar to the
activity of a musician who  exhibits on stages his interpretation of a piece via his
own understanding of the original intention of the composer by reading into the notes.
Probably in the same way as musicians do, mathematicians try to share the joy of
discoveries with the authors of articles, following as closely as possible the lines
given in the works or occasionally  making modifications and others.
At a delivery of a talk, any mathematician is supposed to reconstruct the excitement
of his discovery and try hard to convey it to his audiences. 
\par
This is largely the same with teaching; thus,
a teacher is expected to convey to his students the excitement that he got at the
first moment when he learnt the subject as a student, in addition to a lot of new
aspects that he acquired later as a professional. In fact, teaching a
subject is the best way to learn it; teachers  are probably the most
rewarded in class rooms. Playing a musical piece in front of audiences if small is
said to be far more efficient than practicing it a thousand times  in
solitude. Any layman among audiences can almost immediately judge the quality of any
performance. Just in the same way, students can see how deep or shallow the background
of any teacher is, even if theirs own is awfully poor. I muse on the reason why
such is possible. It is perhaps built in man's very basic faculty.
\medskip
Now, I   turn to the issue of using formulas. It is in fact very
difficult, more than you may imagine, for a mathematician to speak about his own
mathematical interest without brandishing shining mathematical formulas. Mathematics
is such a discipline. Thus, in what follows I shall have to appeal to a little bit
artificial device such as a comparison between mathematics and mythologies. However,
it is also very true that a mathematician cannot conduct his particular  research
solely with formulas, diagrams, and logic, that is, those implements which are
commonly believed suffice to sustain lives of mathematicians.  I need definitely
more, and it is really hard for me to describe what I actually want. I may compare the
situation with that of a musician again; thus, for an instrumentalist it does not
suffice to have the notes, the instrument, and his impressive skill.
\par
Allowing myself to be bold a little bit, it would be more proper for me to claim
that one can do mathematical researches without understanding, for instance, what the 
mathematical logic is. Mathematical articles are armed with logic, yet that 
can never be the aim. In any of my own mathematical discoveries I
was never led to the right destination by logic only. As a matter of fact, there was a
great mathematician who achieved many fantastic discoveries, but yet he did
not seem to be able to wholly comprehend the importance of the logical proof. Thus,
logic is only a kind of excuse, especially for genuine mathematicians, since they
start utilizing it only after having reached what is to be proved. At the very end of
a long arduous journey one finds the truth, suddenly a great view is revealed to
a courageous explorer; the story or the proof is jotted down in a comfortable study
and savoured by armchair travelers. 
\par
It is interesting to discuss, as J.H. Poincar\'e and others tried, how mathematical
discoveries are done, but I think  it is in fact useless, for it does not appear to
me that man is a rational being. Nevertheless, I should add that it is commonly
observed among mathematicians that they use the verb `see' at the very moment when
they sense that the resolution of their problems is imminent. Also the verb `observe'
is often used while conducting researches. A penetrating
observation is indeed a key for resolution. Hence, it is
quite common to find in a mathematical article a grateful attribution of an
observation to a particular mathematician or a group 
even if the observation lacks the logical context. The Buddha, methinks, could observe
all the world instantaneously without any logical deduction; indeed his awful title
is a derivative of the root shared by `observe'.
\smallskip
In passing, I shall relate an interesting story that I learnt at the Tata
Institute of Fundamental Research, Bombay, in 1981 while I was composing my lecture
notes  `Sieve Methods and Prime Number Theory' (Tata IFR -- Springer-Verlag, Berlin
1983). This is an episode of the virtuoso art of observation:
\smallskip
It was then rumored around the world that Russia was about to launch a lunar
rocket to shoot an image of the other side of the moon, which was to be
the first in history. In advance to the event, a famous guru of Calcutta had
been asked to do meditation and draw the picture of the unseen side of the moon; the
result had been shown to a board of scientists including a mathematician, and kept in
a bank vault. When the photograph was actually made and revealed to the public by
Russian scientists, it was compared with guru's picture duly. Amazing! They were
essentially identical. Naturally, the board conferred together and
the guru was asked to reveal the secret. The guru said in reply, ``It is simple. I
extend my spiritual existence to such an extent that the moon comes into
myself. Then I can see readily moon's surface from whichever direction. Your honoured
selves might not believe what I have claimed. However, for instance, a mathematician 
can see a mathematical fact that a layman cannot. A mathematician can,
because he is able to extended his spiritual existence and include the fact into
himself. I do the same but to all directions, while a
mathematician does it to a specific target in a formidably penetrating manner. The
difference is solely there.''
\bigskip
\noindent
{\bf 2. Natural Numbers.} It is said that there are only two kinds
of things in our world. Natural and man-made. You may agree
that the natural numbers $1,2,3,\ldots$ belong to the most natural entities,
as the name suggests well; in fact a great mathematician claimed 
that natural numbers are divine, and
all other numbers are invented. Then, do you know the whereabouts of
natural numbers? They are
like a rigid substance that we may touch; at the same time, they are like the 
legendary aether, for they seem omnipresent in our universe. 
\par
When I muse about this question, I always remember 
a seminal discourse of the Buddha. One day a
man came to him and asked to prove rigorously the existence of the conscience
in all human minds, as the Buddha had preached on it the other day. The Buddha said
in reply, ``Under the full moon, stand with a bowl filled with water. You will
find the moon in it as anybody else does.'' The natural numbers or rather
numbers in general dwell in all human minds, equally and wholly, with no exception and
no incompleteness. Still the whole of numbers dwell in the universe as well.
We are, however, aware of only a certain fraction of such an existence. The human mind
is so deep that even the guru of Calcutta would not be able to fathom it.
In passing, I may note that a bowl seems to have been called `p\=atra'
in   ancient India. It is `ha-ti' in the modern Japanese. Thus,
this pretty common word crossed over the Pamirs, travelled along the Silk
Road eastwards, and came to the minds of my ancestors without much distortion.
\par
By the way, I have been teaching the subject called Complex Calculus for so many
years that I cannot remember when I did it for the first time. At the opening  of
the course in an early week of April, I always start my
lecture with a simple question for my new students. Where are numbers? I of course
expect that someone will give me the right answer. Until several years ago, they had
been perplexed and shown obvious difficulty to reply to my question. In recent years
and this year as well, they gave me the perfect answer, ``In our  minds.'' I suspect
that my questioning has become a ritual. 
\smallskip
Now, I am about to refer to a few excerpts from mythological stories,
for I believe that myths are the very beginning of science or the systematic
knowledge with a firm structure. Ancient people, in the deep past, were already aware
that not only the everyday human matters but also the world of stars and
constellations could not be changeless. For instance, the precession or the
slow gyration of Earth's axis around the pole of the ecliptic seems to have
been found far before the written history of 5000 years. Note that the period of the
gyration is about 26000 years. In any myth, I find a reflection
of the will with which ancient peoples wanted to convey, through and beyond any 
changes in the real world, the eternal truth that they and their
ancestors had extracted from the wealth of observations. In this context, mathematics
is the same, for it is a systematic accumulation of observations of the world
of numbers and a way to convey them to the future. Thus,
the prototypes of mathematical reasonings can be spotted in the construction of
various myths.
\par
A typical example is the definition of man's r\^ole or the creation of man in the
Sumerian mythology. Thus, as I infer from my readings, the gods cannot plough
fields nor harvest crops. However, they have to eat. Thence a principal god
kneaded a block of cray with the divine blood and infused the breath into it to make
man. Man was given the art of agriculture and destined to render thanks to the gods by
enshrining a portion of his crops to the temple. As an obvious consequence, the gods
have to guard man. The point here is to use the r\^ole of farmers, which is visible,
in order to prove the existence of the gods, who are invisible.
The same structure is in the definition of mother's r\^ole that I read long ago in
a myth of an ancient people, whose name I cannot remember unfortunately; 
I surmise that it is of the Judaic origin. Thus, children are
all God's possession; and God is supposed to attend them. However, God is
often too busy with chores and unable to be always beside children. Thence, God
created mother. I am very fond of this definition of mother.
\par
It is never questioned  why the gods need foods nor why children are God's possession.
Just we are left with the great r\^oles of farmers and mothers, and we are
suggested the existence of a splendid harmony that envelopes our lives so gently. I am
awfully impressed by the wisdom of ancient people.
\smallskip
In the modern mathematics as well, it is quite common to push reasonings indirectly.
That is, starting with known facts we assert that such and such exist or be
true, but it is also often the case that any tangible construction of the newly
asserted fact is hard or seemingly impossible, as far as we stay with the argument
employed, though mathematicians are more or less content with indirect reasonings.
There are  many correct statements on natural numbers, for which 
no concrete examples of relevant numbers have been found. Later I shall show a typical
instance with which specialists are frustrated. Anyway, what is essential in the two
myths above is the fact that the reason for that particular choice of the r\^oles of
farmers and mothers is not given at all. This, I think, is a beautiful idea embedded
in the myths. In mathematics, it is precisely an idea or an observation to fix where
the reasoning or the inference should start. The reason why a particular incision has
been applied at the beginning of a discussion  is in the avidya or totally impossible
to explain. True, in any mathematical paper a logical deductions is developed; but
it is often a finely tuned story invented after a discovery.
\smallskip
The mathematical field of researches on integers is called Number Theory. It is quite
hard to make a discovery in number theory. I sometimes compare this situation with
drawing either a revered monk or an apostle. Both are really
strenuous challenges; it is demanded to either investigate or draw anew what have
been done in vastly many manners. Any ordinary
expertise would be defied. One needs genuinely new ideas and infallible 
workmanship together with fierce motivation.
\smallskip
Here is another digression. Discoveries in number theory have certain similarities to
those in physics, a reason for which is in that objects in number theory, i.e.,
integers, are so natural and fundamental that they remind us of atoms 
or more appropriately of elementary particles. Thus, number theory is
sometimes regarded as physics in mathematics; that is, number theory is
thought to be a counterpart in mathematics of physics in natural science. 
I have, however, a certain difficulty to agree fully with this, for numbers that
come up in number theory carry all distinct nobleness, but numbers or constants
in physics do not appear to be so. I am well aware that I am claiming here
something contradicting to what I said in the above. Allegedly, one day a chamberlain
of our last Emperor called a grass at a roadside simply as a weed; he could not
remember the right name of the grass. Immediately the Emperor corrected his poor
subject. ``There is no grass that has the name weed. Each grass has its own proper
name.'' Yes, indeed, exactly by the same token, there is no number which  has no
contents; each has its own assignment in our universe. 
\par
It is strange, however, that we need a great amount of endeavour to establish a
certain number theoretical assertion, in spite of the fact that it must be in
everybody's mind. The whole number system is shared by all people, but yet
sometimes a few hundred years were needed to make a statement rigorous. The late
Prof.\ K. Kodaira once posed me a question. ``Why do we need efforts to understand
works done by others, either new or old?'' He actually compared mathematics with
physics; thus, he asked me also, ``Why can't we see instantaneously
any mathematical fact or discovery like a physical phenomenon?'' At that moment
I could not say any word in reply; if I were asked the same now, I would reply 
boldly, ``Because the human mind is wider than the universe.'' About the depth of our
minds, I shall dwell later. Anyway, I think that
at least number theory is different from physics.
\smallskip
In passing, I should note that number theory is such an old discipline that it
is said to have originated in the book `Arithmetica' that was written by Diophantus of
Alexandria in the middle of the 3rd century AD. It appears to contain a variety of
assertions that are by no means elementary even in the modern standard, though I
haven't had any opportunity to look into the contents.
\bigskip
\noindent
{\bf 3. Primes Numbers.} With this, I now start a story on primes. Thus, the numbers
$$
2,3,5,7, 11, 13, 17, 19, 23, 29, 31, 41,
47, 59, 79, 83, 137, 139, 1481
$$ 
are all primes, since any of these is not divisible by any integer other than $1$ and
itself. The integer $1$ might be seen to satisfy this criterion, but we do not
count it as a prime. The reason for this exclusion will be disclosed later. The
discernment of primes originated probably from the irritation conceived by 
ancient people when they tried to factor a given integer or decompose it into a
product of integers, that is,  execute a perfect division that
leaves no remainder. Indeed, division is the hardest among the four basic arithmetic
operations, as you may find it readily at class rooms of any elementary school. 
The complete  factorization of an integer is obviously harder.
Even today it is difficult to achieve it if a
given integer is huge.  One might argue that the advent of machine
computing has altered the state of affairs. This is partly correct.
However, the change is only in that the size of integers which we may handle
with relative ease has been dramatically increased. In other words, there is a
technological limit with machine computing, and beyond it the situation
is naturally the same as before. We can imagine pretty easily a gigantic number
that any computer is unable to process; obviously it is absurd to
consider the possibility of the converse, for a conventional computer does not imagine
anything.  I may assert safely that the general situation with factorization has not
changed essentially since great antiquity. 
\smallskip
It might sound a tautology, but this difficulty in factoring integers is caused
by the very existence of primes. Thus, if a given integer is a prime, then it is 
of no use to try to factor it. Moreover, it is extremely hard to determine whether a
given large integer is a prime or not. This is called the primality test. Very
recently, there was a remarkable advance in the relevant mathematical field, due to
three Indian computer scientists; however, it does not appear that the difficulty
has been reduced in any practical sense. At all events, the
technological barrier against digital computers will persevere. 
\par
You are to check whether a stone is a ruby or not. You are to spot the exact
location at which a ruby be dug out. Obviously the latter is a far more difficult
task than the former. The primality test is analogous to the former. In what follows
we shall be mainly concerned with the analogue of the latter, that is, we shall
discuss the whereabouts of primes. This is far more difficult than the primality
test, and the field is called the theory of the distribution of primes. 
Rubies and other gems are plenty on the Earth, but yet it is hard to dig them out.
Primes are in everybody's mind, but yet it is hard to locate them
precisely.
\smallskip
I would like to make further discourses before entering into our principal
theme the distribution of primes. Thus, factoring integers seems to have been tried
systematically for the first time in the Sumer-Akkad civilization. Early in this year,
I was at the Archeological Museum of Istanbul, and could see a factor table of
integers written on a clay tablet that is claimed to be some 4000 years old. At the
British Museum and many other institutions are collections of similar clay tablets.
The reason why there have remained such  mathematical documents in plenty is in that
in ancient Mesopotamia  scribes or public clerks seem to have needed mathematical
ability in order to conduct their business, and  the ability to handle division
appears to have been specifically demanded, since division was an advanced arithmetic
operation as is today for most people. Thus, they must have used factor
tables to facilitate their daily works; note that until a few decades ago when 
calculators became popular most engineers had relied heavily on numerical tables of
logarithmic and trigonometric functions. It is not very clear how and in what
circumstances those factor tables were  actually used by scribes. What attracts my
attention is not the technological aspect of those tables but the very existence of
them, that is, the fact that ancient people did know the importance of factoring
integers. As I remarked already, mathematics itself has been rewarded by the
researches on primes; in fact the fundamental discipline Algebra, where the art of
using symbols in mathematical reasoning was first developed,
has an origin in factoring integers. Since it is impossible to imagine any scientific
inference without operations in symbols, I might claim that our science has an origin
in those factor tables of a few thousand years past.
\par
The evidence of such a demand or an encouragement of acquiring the skill of
using factor tables remains well in the curriculums and exercises of the scribal
schools in the Akkadian period. To become a scribe a student had to pass examinations
of mathematics and a classic language, i.e., Sumerian, after several years
of arduous study, as any modern Japanese student has to learn mathematics and English
if he wants to be a civil servant. Mathematics was deified 4000 years or more ago
to be the absolute basics among all educational subjects, 
while the lingua francas have changed along
with the ages and the regional specifications.
\par
One of the most beautiful scenes for us to witness is offered by a child
who does not know yet how to write but is absorbed in inscribing 
clumsy approximations to its name, holding a pencil unnecessarily firmly and 
bending its back. I believe that all
children have fierce wish to learn. This is holy, and must have been the same in
Sumer and Akkad. At ancient schools, boys and girls must have learnt in much the same
atmosphere as at today's schools. At some ruins of wells many clay tablets were
excavated; they were found to contain numerous grammatical and arithmetical errors,
reminding us of modern counterparts.
Perhaps those tablets were thrown into the wells by pupils who did not want to be
further frustrated at their homes. 
\smallskip
Those factor tables are of course
quite incomplete, if viewed with our modern standard. An explanation for this is
that the divisibility by the numbers $2,3,5$ or the first three primes was their
main concern, for other primes did not often come into their daily life, which
is of course the same as even today. Their
society was based on the sexagesimal system, i.e., counting by sixties,  established
already in the fifth millenia BC. The base
$60=2\cdot2\cdot3\cdot5$ is still used in telling time and measuring angle;
note also $360=2\cdot2\cdot2\cdot3\cdot3\cdot5$. Thus, in a
sense, we are as yet in the shadow of the Sumerian civilization, and very
probably will be so in the great distant future as well. I have no evidence to
endorse my assertion but my mathematical experience  suggests strongly that this
choice of the base $60$ was not made by a pure chance. On the contrary, it must have
been done with a far-reaching aim; they made the ever effective decision
based on their mathematical wisdom, employing the first three primes to control
every quantitative aspect of their life. Indeed, for their relatively simple society
the base $60$ must have worked very fine; and locally this is still true with our
modern life. In passing, I add that they used likewise the decimal
system, i.e., counting by ten, but not so widely as the sexagesimal.
\smallskip
The effect of the sexagesimal system or the first three primes are, however, not in
my main concern. I am more attracted by $7$ the next prime. In fact, ancient people in
the East and the West seem to have had a great adherence to the charm of $7$, for
there are so many special words related to this prime, which you may find easily
in dictionaries of any language. For instance, these are pretty well-known: The 7
wonders in various contexts, the 7 hills in many foundation stories of cities, the 7
stars of the Wagon as well as of the Pleiades; also the 7 petals of the lotus
flower on which the Buddha is sitting, and the ceremonies on the 7th day after one's
birth and death in the Japanese tradition. About 280 entries are found in a
Japanese encyclopedic dictionary, all of which begin with the character representing
$7$ and are of certain ancient flavours; more could be found if the words
containing $7$ in the middle are included. 
Why is this preponderance of stories and words  with the
prime $7$? As to the next five primes, in the same dictionary, there are 9
entries with the prime $11$; $26$ with $13$; $9$ with $17$; $7$ with $19$; $1$ with
$23$. One striking example is in the Esoteric Buddhism;  there is an instance in
which the three primes $7$, $17$, and $23$ are closely related.
\par
Certainly the grandest and probably the
oldest among those myths involving $7$ is the Deluge Myth in the Epic of Gilgamesh. 
According to the the standard Babylonian version of the epic 
(transl.\ by A. George, Penguin Books, London $2003$), the great story was
revealed to Gilgamesh the king of Uruk by Utanapishtim the immortal human living in
the far distance. The tempest for 7 days on the Sea and the arduous 7 days on the
Mount Nimush, with Utanapishtim, all his kith and kin, members of every skill and
craft, and the seed of all living creatures on board the Ark.
Moreover, the denial of the immortality for Gilgamesh,
failing to prevail against sleep for 7 days.  
\par
The discovery of the Deluge Myth among the tablets from the great Nineveh
library of king Ashurbanipal was done by G. Smith in 1872 when he was inspecting the
eleventh tablet of the Epic of Gilgamesh, at the British Museum. It is recorded that
at that moment he was so greatly excited that he started removing his clothes. He had
been an apprentice of a bank note engraver, but learnt to decipher the cuneiform
scripts by himself. His talent and
enthusiasm had let him win a position, and he left the great
achievement. His is indeed a beautiful life. Only the passion for learning could
make it possible.
\smallskip
Then, I ask where primes had their genesis. It is marginal to determine when
the word `prime number' was used for the first time. The essential is to ask in which
civilization people cherished and played with primes for the first time in man's
history. The reason why I think this is a major issue was already suggested
above; the answer should locate where the modern science originated.
Expressing this in a somewhat different way, I think that the tradition of science
supported by abstract reasoning must have started with segregation of certain primes,
if only a few, as basic constituent elements of the whole number system. Similar
segregation was certainly applied to natural substances as well in ancient
civilizations. However, the isolation of certain integers as something special is far
more a mental faculty, and the result is definitely universal, for this is
purely the fact in everybody's mind. Why has the sexagesimal system
survived till today? Because it was the result of an abstract thought targeting
a grand general dwelling in everybody's mind. 
\par
It is true that the concept of primes as a whole is due to ancient Greeks, and is
of course tremendously important, as I shall dwell later.  However, I do not think
that devising a concept is of absolute importance, for the birth of a concept
needs a mother, and that is a poignant special, pregnant with important consequences.
Thus, who played with prime numbers for the
first time? Who enjoyed the first moment of the true learning or man's
greatest luxury? Who realized for the first time any existence that is special as
well as universal? Since my boyhood until recently, this query about primes had
charmed me. I am about to relate how it was resolved.
\smallskip
It was around the end of October last year, when I visited the Maruzen bookstore
newly opened close to the Tokyo Central Station. I came to the section
of mathematics, and found K. Muroi's book `An Introduction to Babylonian Mathematics'
(Univ.\ Tokyo  Press, Tokyo 2000 (Japanese)).
I picked it up, and opened the index as I do always when I want
to see the quality of any academic book. I found `prime' among the entries. 
I was, however, not much interested in the fact that the word was there, for
I believed already that the concept of primes must have been acquired by the ancient
Greeks from the greatly precedent civilizations in the Near East. Thus I returned
the book to the shelf without looking into the relevant part, and went back to
my home. This might be a useless addition but
I have the impression that most popular books on primes, either in Japanese or in
English, are unfortunately not very satisfactory, although I understand how difficult
it is for those authors without any experience of researches on primes to describe
this  highly specialized subject in plain words. Because of this, I did not
see the contents of Muroi's book, although I was quite aware 
that it was never for novices. However, already on the train of the Yamanote Line
that I took returning home, I became increasingly attracted by the book, and further
more at home. In the next morning, I came to the bookstore again, and scanned
the section on primes. I really felt the secretion of the adrenaline into
my veins. Immediately I bought a copy, went back to my home, and enjoyed a
nightlong reading.
\par
Muroi excavated a tremendous treasure from the 4000 years accumulation of 
the sand called time. Mathematicians in ancient Mesopotamia are known to have
treated various equations with ease. What is amazing in Muroi's discovery
is the fact that primes were found concealed in the
coefficients of equations as their factors; those primes are listed
above. In a computation of a square root, the prime $1481$ is embedded. Can you
check swiftly if this is a prime? It is apparent that primes were used to make
equations look more interesting or mysterious. The intriguing device by
ancient mathematicians penetrated my mind at the middle. They are surely among the
first people who genuinely played with primes.
They did mathematics with primes, quite in a modern sense. Therefore, my quest for the
genesis of primes ended with Muroi's book. Resplendent gems are inlaid in mathematical
equations! I think that Muroi made a discovery greater than H. Carter's at the tomb
of King Tutankhamen, for those gold and jewel treasures are unquestionably great, 
proudly exhibiting the unbelievable  artistic and technological achievements of the
ancient civilization, yet they are perishable in the great extent of time. The
primes inlaid in equations will never perish. Indeed, the mischievous, witty spirit
of those mathematicians who were active 4000 years or more ago will be fresh for
ever, endorsing magnificently that pronouncement in the preface 
of the Kokin-waka-shu. 
\smallskip
At the end of his book Muroi cast a sympathetic glance toward the citizens in ancient
Mesopotamia where peace was rare and always short-lived. 
Perhaps because of that, people made trips to mathematics in the quest of peace and
eternity. This reminds me of the life of my late mentor Prof.\ P. Tur\'an a
great number theorist. Once he told me about his life that was full of turmoils,
caused by wars and enmities. I asked him how he could survive the hardship. He smiled
lightly, and said, ``Doing mathematics. We had seminars at the institute regularly,
even during the uprising against the Soviets.  To be absorbed in mathematics was the
sole remedy for me.'' It was early 1970's, and I never imagined that I would hear
about this invaluable virtue of mathematics again in 1999. Then my great friend A.
Ivi\'c was in Belgrade under the air raids. We could maintain the e-mail connection
during the war, though it was not very reliable because of frequent power cuts over
there. One day he sent me such a fantastic result  in a {\TeX} file that I could not
readily believe. In fact I could not believe the fact that he could do mathematics
under that hardship; until the early spring of that year, it had been absolutely
impossible by any stretch of my imagination that any of my friends was being involved
in a war. Later Ivi\'c told me that doing mathematics had been a remedy for him
during the war.
\smallskip
In passing, I may stress the virtue of libraries as well. Those ancient
royal libraries and modern digital archives are all the same in their
sublime purpose, that is, to store anything written in letters, symbols, diagrams,
images, etc.,\ and to pass them to the future generation without any failure.
 Ashurbanipal, who was a rare literate among Assyrian kings, must have
thought that future people would be grateful to him for his great project of a
complete library. Indeed, his effort was appreciated 2500 years later. I wonder how
the situation will be in the future of 2500 years.
\smallskip
In the above I mentioned frequently ancient Mesopotamia, which I am now leaving
for ancient Greece. The reason why I am fond of thinking about the civilization 
of the land between the great rivers is in that I can see well people clinging to a
small peace in small corners of bustling cities. That is quite similar to my
life in Tokyo.
\bigskip
\noindent
{\bf 4. Sieve Methods.} Thus, I have come to ancient Greece, perhaps from Levant by a
Phoenician ship; and further to the metropolis Alexandria the
greatest academic crucible of antiquity.
\smallskip
At the beginning of all
modern textbooks on number theory, you should find a theorem asserting the unique
prime factorization of any integer. This means that any integer
can be expressed as a product of primes, e.g., $546=2\cdot3\cdot7\cdot13$, and if you
find such an expression by any reasoning, then you do not need to look for any other.
Of course the order of prime factors is disregarded, for instance we could write
instead $546=13\cdot2\cdot7\cdot3$. 
In the `Elements' which was compiled by Euclid of Alexandria around 300
BC based on various older sources, a theorem is finely proved, from which this
uniqueness assertion can readily be deduced. Moreover, it is clearly established that 
there are infinitely many primes; yes, infinitely many gems in your mind.
It appears, however, that the uniqueness of prime factorization itself is not stated
anywhere in the book; for this observation, I thank K. Saito a mathematics historian.
At any event, it is historic that with Euclid's book the whole primes entered
into man's vista, and the concept of prime numbers was firmly established. 
\smallskip
My main purpose in this article is not to discuss the history of mathematics; and thus
I shall not disclose my thoughts about the r\^oles of ancient Babylonians,
Greeks, and other peoples concerning the concept of primes. Nevertheless, I should
express my difficulty to accept the claim, still popular among non-specialists,
that the ancient Greek mathematics burst out of the blue. At least it is out of the
question to assert that the seed was solely in ancient Greece. However, the logic in
the Greek mathematics is indeed beautiful and grand, reminding me precisely of Greek
temples and sculptures. As I wrote above, I was at the Archeological Museum of
Istanbul early in this year. There I was really impressed by the perfect beauty of
the relief on the sarcophagus of Alexander the Great, which is actually attributed to
a king of Sidon, though. The origin of the word sarcophagus and the purpose
of the item are gruesome, but I was absorbed, for I had never encountered
anything of such a perfection in any of my somewhat extensive experiences at many
museums. It was made more than 2300 years ago; the sculptor was essentially a
contemporary of Euclid.
\par
It should be stressed here that the original volumes of Euclid have not survived 
as Greek sculptures; thus,
the modern discussion about the achievements by Euclid and
his school is based on editions and translations compiled far later, well
1500 years later at least. On the other hand, those mathematical clay tablets from
Assyria are genuinely original, that is, no editorial works were done on them.
One should keep this great difference in his mind when attempting any comparative
study of ancient mathematical traditions. 
\smallskip
Whichever ancient civilization he
prefers, anyone must have difficulty to live in archaic styles. The Babylonian
mathematics is highly seminal but too archaic; and the
grand Greek mathematics is the same. If you think a
little deeper about the two Greek assertions in the above, you should feel much
anxiety. Well, there are infinitely many primes; then, is it possible to say how much
they exist in plenty, in the practical sense? That is, even if this might be useless,
is it possible to count them one by one, to assure ourselves that they are indeed
plenty? Also, the expression of any integer as a product of primes is always possible
and the mode is unique; then, how can one find such a factorization, actually?
Unfortunately, the logic employed by Euclid does not supply any hint to resolve these
two basic problems.
\par
I said boldly in the above that a logical proof could be an excuse. Indeed,
mathematicians often realize only after publishing articles that they have not
understood fully what they asserted and proved proudly in their works. To establish
an assertion with mathematical logic does not imply immediately that the statement has
been fully understood. A proof guarantees only the correctness of the relevant
assertion. That might sound bizarre, but I myself have such experiences with my
publications. In much the same token, the author of the Elements looks similar to the
majority of modern mathematicians, for he touched only the surface
of the matter. Yes, I know I am too harsh, for those two basic problems on primes are
still lingering nightmares in mathematics. Perhaps, it is more appropriate to praise
Euclid for his achievement, given the mathematical background of his time. Successors
must be modest.
\medskip
My aim in what follows is to relate man's struggles against the
enormous difficulty in counting primes. This theme has, however, the tradition of 2200
years at least, and I have to restrict myself to really pivotal events,
in fact, in the last 150 years mainly. 
\medskip
Probably you heard about the Sieve of Eratosthenes some years after
graduating from elementary school. You leave the integer $2$ intact, and erase those
integers divisible by $2$. The integer $3$ remains. Leave it, and erase those
remaining integers which are divisible by $3$. The integer $5$
remains. Leave it,  and erase those remaining  integers which are 
divisible by $5$. The integer $7$ remains. Repeat the same with $7$, $11$, $13$,
$17$, $19$, $\dots$, so on. This procedure is attributed to Eratosthenes of
Alexandria, active in the later half of 3rd century BC. It reminds us of sifting;
thus the above term. 
\par
Now, when the sifting has ended, all integers remaining are primes, and all
primes appear there duly. Namely, the table of primes can be constructed
by this procedure. The logical confirmation is easy. It appears that Eratosthenes'
motivation was to complete those incomplete factor tables that I mentioned above, or
rather he wanted to construct the table of the prime factorization. The construction
of the table of primes does not seem to have been his intention; and that  aspect of
his method was stressed only 1500 years later. Of course they are the same, since,
marking integers appropriately  instead of erasing, the complete factor table  becomes
visible. For instance, $1092$ is marked with $2$ for the first time. Since the
quotient $546$ is smaller, it may be marked fully already; if not, apply the procedure
to it. In this way we can achieve a full decomposition or marking of $1092$. The
integer $1481$ will get no mark, and it is a prime.
\par
Despite this common attribution to Eratosthenes, it appears likely 
that the procedure had come from the Orient, in view of the long history of
factor tables and on noting that it can be conducted independently of
Euclid's assertions on primes. Ancient mathematicians made academic trips worldwide
as  their modern colleagues do today;  an Ebla tablet tells such an
episode that took place more than 2000 years earlier than Eratosthenes' time. 
Eratosthenes was a poet, astronomer, and director of the great library of Alexandria.
He is said to have been old and unable to enjoy books when he died on voluntary
starvation. 
\smallskip
Now, as I promised above, I shall reveal the reason why the integer $1$ is not
counted as a prime. This is another digression, and a somewhat
old story from my experience. I was then going abroad. When I came to my seat on a
plane, it was soaked. Embarrassed, I indicated  my trouble to a stewardess. A little
later she came back, and offered me a seat in the first class. Wonderful. I felt
slightly uncomfortable with my sudden luck, and to calm myself down I started
checking my article for my talk at a conference. After a while, the gentleman
in the seat next to mine spoke to me in a heavy German accent, ``I am a physicist. Are
you a mathematician?'' ``Yes, I am.'' ``What is your specialization?''``I am
specialized in the theory of the distribution of prime numbers.''``Hum, is there
still any unsolved problem about prime numbers?''`` Yes, a lot,
indeed.''``Well, there is one question that has been lingering on in my mind since my
boyhood.''``What is that?''``Why isn't the number $1$ a prime number?''``Do you know
the procedure called  the sieve method?'' ``Yes, I think I learnt that at Gymnasium.
It starts with erasing integers divisible by $2$. Right?''``Fine. Then, what will
happen if you begin the procedure with the number $1$?'' After a thoughtful pause, he
exclaimed, ``Oh,  Danke!''
\smallskip
The correct explanation is that if we count $1$ as a prime, then Euclid's 
assertion of the uniqueness of the prime factorization would fail to hold. However,
this point is not very relevant to our present purpose, and let me return to a more
serious issue that I touched on already. Thus, let me ask you whether you find
it strange or not that in spite of the existence of the table of primes the research
on the distribution of primes is still going on. A table is given to us;
then our common sense tells us that the distribution of the entries on it should be
obvious.  In the case of the table of primes, this does not hold unfortunately,
at least for the present. A somewhat
naive explanation of this defiance to our common sense is that the method of
Eratosthenes is not practical, or expressing the same a little more precisely, 
it does not allow us to foretell a priori where a prime appears. 
Eratosthenes' sieve yields the table but yet it
cannot be used to pinpoint each incident of the entry. This sounds too
bizarre, however.
\par
Here the greatest theme of mathematics looms over us. That is, we are
facing the riddle of the infinitude. It is known that the above sifting is effective,
if we are given a certain finite range where the procedure is to be applied;
however, if the range is indefinite, then Eratosthenes' sieve becomes almost
powerless. In the above, I wrote ``Repeat the same .... Now, when the sifting
has ended, ...'' There I devised a kind of deception. 
An infinite repetition of dividing infinitely many times must have been
accomplished when the sifting has ended. Is this possible? 
Yes, but only in our minds. It is impossible in reality. Namely, it is only
via our mental faculty that we can see that the sifting proceeds 
and the table of primes grows,  beyond any prescribed finite boundary. The table is
ever incomplete or always in the making but once you impose a finite boundary 
it can be readily or has been already completed up to that point. The complete table
exists only in the process of `mathematical  limit', which takes place in our
fathomless minds, nowhere else.  Discussing natural numbers, I remarked that often in
mathematics  ``any tangible construction of the newly asserted fact is hard or
seemingly  impossible.'' The construction of the table of primes is a typical
instance, as far as we appeal to Eratosthenes' sieve. A complete practical table is
yet a myth, although I should not deny its possibility in the future, well, in the
greatly distant future when one may dig out primes one by one as pebbles, beyond any
boundary.
\smallskip
Probably because of this ineffectiveness, Eratosthenes' sieve had been disregarded
for more than $2100$ years, until V. Brun  revolutionized number theory in
$1915\!\sim\!20$, quietly with his new way of sifting. As to the counting of primes,
there had been a great advance two decades earlier; but that will be related later.
Brun discerned that Eratosthenes' sieve was too complete to be effective or
practical.  His idea was to make Eratosthenes' sieve flexible by introducing certain
skips in those repetitions of sifting; namely, the erasing or the marking is not 
performed with all primes you encounter in the procedure described above. This must
sound strange, because it is not clear how you are to pick up primes to
use in Brun's way of sifting. I should have made
it precise that in Eratosthenes' sieve you get a new prime automatically at the
beginning of each repetition of sifting, but Brun assumes instead that you are given
a finite number of primes prior to all the repetitions of sifting. I am about to
provide you with more details. You
need to be prepared for a dose of a small amount of mathematical reasoning.
\par
Thus, if you study the procedure above a little
more,  you will find that when the round with the prime $7$ starts, you have
processed the range up to 49 already; that is, all primes less than 49 have been
found.  To see this, let us pick up an integer less than $49=7\cdot7$. If it is  a
prime, then there is nothing to discuss further. If not, by definition it can be
divisible by an integer larger than $1$ and less than itself; of course the
r\^oles of the divisor and the quotient are exchangeable, and 
we see that the integer under question is in fact expressible as a product of two
factors both larger than $1$. Then, one of these two is less than $7$; otherwise, both
are larger or equal to $7$, and the product must be larger than or equal to $49$,
which contradicts to the initial assumption that we are working with integers 
less than $49$. Hence the integer under question is divisible by an integer larger
than $1$ but less than $7$, that is, it must have been already 
sifted out with a prime less than $7$, i.e., one of
$2,3,5$. In other words, there is actually none to erase less than $49$ at the round
with the prime
$7$. This ends the discussion of the case with $7$. By just the same reasoning, when
you start the round with
$37$, which is a prime, you have already finished the sifting up to
$37\cdot37=1369$. With the prime $1481$, the sifted range must be up to $2193361$
which is large indeed. 
\par
Therefore, you need primes less than the square root of the
upper bound of the range you want to sift, which is a little theorem. 
To list up all primes less than $100000000$, you need all primes up to $10000$, and to
have the latter you need to know primes up to $100$, which is easy. With this,
probably you will doubt the truth of my claim that Eratosthenes'
sieve is ineffective. One way to gain your consent is to estimate the number of
operations actually involved, but I leave the task to you. You will certainly be
convinced of my claim. Thus, I may assert that Eratosthenes had in fact so many primes
in advance, i.e., all primes less than the square root of the prescribed bound, that
his sieve could be a method for the primality test within the bound,
though it involves a too huge number of applications of  division to be
regarded as a practical method.  
\par 
Brun gave it up to have so many primes in his hand as
Eratosthenes did prior to sifting. You will
anticipate immediately that his device should be destined to fail in the exact
primality testing. Yes, it is; and the result of his argument is
rendered  with inequalities instead of the equalities which is the case with
Eratosthenes' sieve. He could estimate, however, the number of those integers in a
prescribed boundary that survive his application of sifting; or more precisely, he
could  choose those primes in his hand in a highly effective way so that  his
estimation became remarkable, the more the range grew wider. With Eratosthenes' sieve
such an estimation  becomes meaningless as the range grows indefinitely, as I have
indicated above.
\smallskip
Brun was the first who found how to hammer pitons into the formidable bluff of the
infinitude of integers. His way of sifting is termed Brun's Sieve. It appears that
his idea was too bizarre for his contemporaries to swallow. Nevertheless, it was
appreciated widely decades later, and brought eventually a Sturm und Drang to 
number theory, especially to the theory of the distribution of primes
that lasted for more than 60 years. Of course the movement is still going on. 
\par
The same as most breakthroughs in mathematics, his idea bears a primordial 
simplicity. Simplicity is a power. Eratosthenes could  have devised the same.
When I encountered  Brun's sieve for the first time, I was reminded of the legend
that Alexander the Great cut with his sword the intricate knot of Phrygian King
Gordius, and proceeded to Asia. In fact, in my mind Brun is mightier
than the great king, for he cut the enigmatic knot that had survived 2100 years
without any sign of wear. I am aware that I might cause misapprehension but I dare
to note that Brun's idea shares a basic philosophy with the revolutions in 
other fields of mathematics and in physics, all of which took place in the early 20th
century, that rigorism be overridden to attain flexibility, even if ambiguity
or uncertainty is brought in. 
\par
I wonder why the classical rigorism or the
kingdom of equalities had to give
way to methods tinted with statistical reasonings. 
Once M. Jutila, another great friend of mine, rendered
this situation with a quite smart maxim of his. Thus, equality might contain
noise more than inequality does. Applying this to sieve methods, I can say that
Eratosthenes' sieve contains noise because of its exactness, and Brun devised
a filter to segregate the essentials while discarding some part of noise out. Of
course the filter should be chosen in an optimal way. Brun's choice of primes
in his hand was indeed a remarkably effective filter, being neither too tight
nor too loose.
\par
As another digression, I may point it out that at the basis of the modern calculus,
i.e., the art of differentiation and integration, is an extremely important
reasoning consisting of a combination of inequalities. Thus, the mathematics 
taught at university demands students be familiar with inequalities; indeed, 
they are to be treated as a constituent of general citizenry who have to 
read inequalities more often than equalities, not only in
mathematics.
\smallskip
The above question about the great scientific changes 
will look in fact quite relevant to our present discussion, if I dwell the 
developments after Brun. Thus, the statistical reasoning has got a tighter
grip on the theory of the distribution of primes; 
key words are, for instance, the large sieve of Yu.V. Linnik 
and the $\Lambda^2$-sieve of A. Selberg. However,
I should not stray from the main purpose of mine. You are referred instead to my
article `An overview of sieve method  and its history' [arXiv: math.NT/0505521] to
appear in  Sugaku Expositions AMS, in which  you will find the modern situation of
sieve methods with considerable technical details. 
\smallskip
Nevertheless, I should touch on the Twin Prime
Conjecture and the Goldbach Conjecture, which are both still open, and too famous to
be left untouched by any account on primes.  Twin primes are two primes which are 
different by $2$; for instance $41$ and $43$ are twins in this sense;
likewise $617$ and $619$ are. The
conjecture asserts that there are infinitely many twin primes. The Goldbach Conjecture
asserts that all even integers, i.e., integers divisible by 2, which are  larger than
$2$, can always be expressed as a sum of two primes as $100=47+53$ and $1000=281+719$.
Both belong to the greatest problems in mathematics. You must be
amazed by the fact that problems of this simplicity have defied tremendous efforts of
best  mathematicians for some $250$ years. It is in fact pretty easy to formulate
problems involving primes that may intrigue specialists and cause insomnia as these
conjectures or perhaps more deeply. One could regard this as indicating mathematics
is still pretty immature. Well, it might be the case, but my view is different. I
think that the situation implies simply that we are facing a terra incognita lying in
anybody's mind. Yes, again, man's mind is fathomless.
\par
Brun's original motivation was to modify
Eratosthenes' sieve so as to close in on these two great problems; and indeed he
achieved a remarkable progress. He proved that there are infinitely many pairs
of integers different by $2$ both of which have $9$ prime factors at most.
This was then an astounding approximation to the Twin Prime Conjecture. He proved
likewise an analogous assertion for the Goldbach Conjecture. It should be stressed
that Eratosthenes' sieve is unable to yield anything similar.
\par
After Brun's achievement, the progress was slow but it was along the line
set by him, and a climax came in $1966$, when J.-r.\ Chen announced that
he had proved that $2N=p+P_2$, namely any large even integer $2N$ can be
expressed as a sum of a prime $p$ and an integer $P_2$ that has two
prime factors at most. His argument yields also that there are 
infinitely many primes $p$ such that $p+2=P_2$. These are tremendous progress
towards the resolution of the two great conjectures; perhaps, just one step
short of the final goal. 
\par
I learnt Chen's announcement quite early.
It was in winter 1966; I had already moved to Tokyo from Kyoto, thanks to a
position offered by Nihon University, where I am still teaching. There was
a colloquium talk by S. Uchiyama at Tokyo University. After his
talk, I had a short discussion with him. He said, 
`An astounding announcement has
been made in China.' `What is that?' 
`A mathematician named Chen Jing-run
has claimed $p+P_2$ for the Goldbach Conjecture.' `Anything about
his method?' `Mean Prime Number Theorem of E. Bombieri and I.M. Vinogradov, and two
lemmas, but I have difficulties with one of the latter.' Then he gave me
a copy of the now famous announcement.  
Actually, the two-page announcement was incomprehensible for me who
had no experience with any subjects Uchiyama had kindly mentioned. In fact, that
was the moment when I decided to study  sieve methods; but all
essentials were dispersed among isolated articles in a few languages, which forced
me to learn Russian in particular. 
\par
A turmoil caused by a deplorable 
political crime, termed  inconsistently the Cultural Revolution
($1966\!\sim\!1976$), had already been spreading, and
Chen would lose seven years until the publication of his full proof. 
In a year, the Student Revolt started throughout Japan as in the
western Europe; my university and the Surugadai district were the stage of
battles between students and the police for nearly a year. I could not do
any research; I felt ambivalence toward uprising students. Another year later,
I was on PS Baykal leaving Yokohama for Nakhotka. Tur\'an had warmly accepted my
plea to study the distribution of primes under his supervision; I was
in part fleeing to Budapest. It was in January 1970. Pack ice hitting the ship; 
in a freezing night of Khabarovsk, young
women in rugs working at railway with heavy tools; their stern
eyes toward me; Moskva in a record cold spell; the barren Red Square and the golden
Kremlin; the Beograd Express via Kiev; tea from samovar served by a young woman
conductor; soft Russian melodies from somewhere of the car; the frozen window of the
compartment, red with the dazzling icy sunset; lifting wagons and changing bogies
at Chop the border to adjust trains to the gauge difference; the Keleti station
of Budapest, so dark as blackout. I arrived my second town, dispirited. However, in
the next morning I was absorbed by the great view of the white V\'ar and the frozen
Donau with splendid bridges. 
\par
Since then, I have made
numerous trips to meet with mathematicians and experiences. Naturally, the first
journey to Budapest is the most memorable. Likewise unforgettable are the stay at the
Tata IFR and the excursions to the south and the north India. The chill of the stone
bed in a chamber at a rock monastery of Ajanta pleased me, as Buddhist disciplinants
must have felt 2000 years ago; the place was still filled with the spirit of
learners. On the other hand, descending to the earth, among the dishes I have ever
enjoyed on my trips, the best is the quince compote at a restaurant on the shore 
of the Black Sea where the Bosporus starts south.
\bigskip
\noindent
{\bf 5. Riemann Hypothesis.} Now I shall relate the very basics from the theory
of the distribution of primes. I shall have to rely on mathematical
formulas and notations more than before. If you have
no mathematical background including complex calculus, then this section should be
hard for you to fully comprehend. Nevertheless, I shall try to make my discussion as
plain as  possible so that you may try to follow me.
\smallskip
I shall first need the notion of the complex numbers. It is
usually introduced to students when they learn the general way to solve quadratic
equations. However, you are not required to know it precisely; the essential is to
keep in your mind that there is a wider world of numbers including the real
numbers, i.e., those you use daily. It is always possible to compare two real
numbers; that is, one comes either before or after the other. Perhaps you agree that
because of this the real numbers seem to  dwell in a straight line. We adopt this
common sense as a mathematical convention; thus the world of real numbers is
identified with a particular straight line, termed the real
line, which is in everybody's mind. On the other hand, the world of complex
numbers is identified with a particular  plane that spreads including the real line;
any plane is also in everybody's mind, and thus complex numbers have long dwelt in
your minds already. This plane is called the complex plane. 
Between two points on a plane you have difficulty to say which comes after the other,
unlike on a straight line. However, on a plane you can instead move freely, while on a
straight line your movement is strongly restricted. Thus, in the world of complex
numbers or on the complex plane you may enjoy a great freedom. This freedom in
mathematics was discerned early in 19th century by most European mathematicians, and
it was a start of the modern mathematics. One of the aims of the mathematics taught at
university is to enjoy this freedom as much as possible; and the most
relevant subject is called complex calculus, which I have mentioned already that I
have been teaching for many years. There is a branch in number theory that depends
largely on complex calculus. It is called  Analytic Number Theory, to which I am
specialized. This is a mathematical tradition established by G.F.B. Riemann
($1826\!\sim\!1866$), as I shall indicate more later.
\medskip
With this, let us denote by $\pi(x)$ the number of primes
less than a given real number $x>0$. It should perhaps be stressed that this $\pi$ has
nothing to do with the well-known number $3.141592\ldots$. Our $\pi$  stands
just for our particular function, that is, $\pi(x)$ represents the way how a moving
real number $x$ or a real variable determines a particular value; for instance,
$\pi(20)=8$ and $\pi(10^8)=5761455$, where $10^8$ denotes the $8$-fold multiplication
of $10$ or a hundred million. There is no difficulty to fix the former, but the latter
demands a huge number of steps as I have suggested already; the value $5761455$ is
obtained by actually counting all primes up to
$x=10^8$ one by one. That is too tedious, and naturally one may wonder if there
exists a practical way to compute $\pi(x)$ for any given $x$. This is the most
fundamental problem in the theory of the distribution of primes. It is true that
there exists a complete table of primes in our minds, but it is only a mathematical
limit; that is, only its existence is certain. Under a situation as this,
mathematicians ponder usually on the possibility of devising a way to compute
approximative values. Namely, they try to seek for a method that yields the
value of $\pi(x)$ for any $x$ within a certain admissible error. It is
indeed amazing that such a method does exist. It is via the logarithmic
integral, denoted by $\r{li}(x)$. You do not need to know the rigorous
definition, but anyway it is given by the integral
$$
\r{li}(x)=\int_2^x{du\over\log u}\,.
$$
What is essential with this is in that there is no need to know anything about
primes in order to compute the value of $\r{li}(x)$ for any $x$; moreover, there
exists a highly practical way to compute it. For instance,
$\r{li}(10^8)\fallingdotseq 5762209.375$ is fairy easy to obtain. This is close
to $\pi(10^8)$; the relative error is $0.0131\%$, i.e., $|\pi(10^8)-\r{li}(10^8)|
/\pi(10^8)=0.000131\ldots$, which is almost negligible against our daily standard.
Here the notation $|\ast|$ means that we are concerned solely with the size
of $\ast$ a certain quantity either positive or negative; absolute value is
the term. It was K.F. Gauss who first observed
that $\r{li}(x)$ might be close to $\pi(x)$ for all $x$; he was then just a boy.
Naturally he checked only those $x$ much smaller than $10^8$.  
\smallskip
Then, is it really the case that
$\r{li}(x)$ is an approximation for all $x$? That is, is it true that the error
$|\pi(x)-\r{li}(x)|$ is ignorable compared with $\pi(x)$? We should of course make it
precise in what sense it is ignorable. Thus, let us check
how the relative error varies along with $x$. It is 
about $0.00334\%$ with $x=10^9$, and $0.0007\%$ with $x=10^{10}$ or ten billion.
Most probably you will surmise that the relative error should tend to naught.
In fact it had remained as Gauss' Conjecture on $\pi(x)$ until
$1896$ when it was proved
that in fact the relative error tends to naught. It was done independently by
J. Hadamard and  Ch.\ de la Vall\'ee-Poussin. Close to the end of the 19th
century, man finally reached a really tangible result in the art of counting
primes. This monument in mathematics is rendered in the statement
$$
\hbox{The Prime Number Theorem:}\qquad \pi(x)\sim\r{li}(x)\quad 
\hbox{as\quad $x\to\infty$},
$$
where the symbol $\sim$ stands for the fact that the relative error tends
to naught.
\smallskip
You may argue that it is more important to know how this was made possible. In fact,
all basics that Hadamard and de la Vall\'ee Poussin needed had been prepared
by Riemann in his paper `\"Uber die Anzahl der Primzahlen unter einer gegebenen
Gr\"osse', i.e., On the number of primes less than a given bound (Monatsber.\
Akad.\ Berlin, 671--680 (1859)), which is certainly one of the greatest mathematical
works of all time, despite its unimpressive size and its awful lack of
details. You see here that the completeness of the mathematical logic
is not always necessary for mathematicians to attain respect, as I said before in
other context.  The core of Riemann's idea was in his introduction of the
zeta-function. The same had in fact been employed by L. Euler and P.L.
Chebyshev in their researches on $\pi(x)$; but because of the  incomparable
importance of Riemann's contribution, it is now termed the Riemann zeta-function, 
and denoted by $\zeta(s)$. Unlike his predecessors, Riemann considered 
$\zeta(s)$ with a complex variable $s$, that is, to each value of
$s$ moving around on the complex plane a complex number $\zeta(s)$
corresponds. The correct definition is simple for anyone who learnt
mathematics at university, as it is given by
$$
\zeta(s)=\sum_{n=1}^\infty {1\over n^s},\quad \Re s>1.
$$
The uniqueness of the prime factorization of integers implies that
$$
\zeta(s)=\prod_p\left(1-{1\over p^s}\right)^{-1},\quad\Re s>1,
$$
where $p$ runs over all primes, and $\prod$ means product as $\sum$ stands for sum.
Here,  I should stress again that you do not need to know what actuary these
expressions are; thus the condition `$\Re s>1$' is left unexplained. What you need to
observe is that the latter expression, which is originally due to Euler and thus
called the Euler product, exhibits that $\zeta(s)$ is directly, beautifully
related to the whole of primes. Hence, it may be natural
to expect that the property, which professionals call the analytic nature, of
$\zeta(s)$ as a function of a complex variable should be related to the distribution
of primes. In this context, what is further remarkable is in the fact the $\zeta(s)$
has the first expression. Because of it, we can consider
$\zeta(s)$ throughout the complex plane where $s$ moves around, while the sole
use of the Euler product cannot  yield the same. The latter point might be bizarre.
However, we have to accept this as a fact or just that our $\zeta(s)$ is many sided
as usual with most of existences in our universe, and each expression is only a part
of the whole. Indeed, it is a great luck that the first expression never suggests any
close relation of $\zeta(s)$ with primes, whence we can analyze $\zeta(s)$ without
referring to primes. Anyway, from now on
$\zeta(s)$ is defined for  any complex number $s$; this is a typical instance of the
procedure called analytic continuation in complex calculus.  
\smallskip
Now, the great discovery of Riemann is in the diagram
$$
\pi(x)\longleftrightarrow \left\{\hbox{complex zeros of 
$\zeta(s)$}\right\}.\eqno(\alpha)
$$
Here, that $s=\varrho$ is a  complex zero means that $\zeta(\varrho)=0$, and
$\varrho$ is not a real number, i.e., not on the real line. Again you do not
need to understand this diagram very precisely, though of course I would
appreciate if you do. You are only required to have the image in your mind that
counting primes is closely related to a special set of complex numbers which
are called complex zeros of the Riemann zeta-function. This fact is indeed the
duality in the research on primes that I touched on at the beginning of the
present talk. Perhaps, $\zeta(\varrho)=0$ is the moment when $s$ has walked around
looking for primes and the poem $\varrho$ has just scintillated; well, I am pretty
aware of a forced analogy. Nevertheless, this description of
$(\alpha)$ is rather close to the mathematical reality, for the complex zeros are
defined through the first definition of $\zeta(s)$, that is, without touching 
primes. 
\par
You may wonder what merit is in this translation; in fact, $(\alpha)$ appears to be
only a replacement of $\pi(x)$ by a set of more difficult, mysterious numbers. You
are expected here to remember that complex numbers dwell in a plane. Yes, the right
side of
$(\alpha)$ belongs to the liberal camp; and then you may suppose that more implements
should be  available in dealing with the flexible right side than the traditional
left side which belongs to the archaic conservative camp. This is indeed the case, as
Hadamard and de la Vall\'ee Poussin demonstrated. 
It should be added, however, that about half a century
later Selberg and P. Erd\"os achieved 
an elementary proof of the prime number theorem, that is, without
flying into the world of complex numbers. Their proof appears
hard and inflexible  to me, although I certainly appreciate their
achievement. By the way, `elementary' is a misleading term, for it 
does not mean `easy' in our context. It is used when one proves any assertion on
integers without enjoying the freedom on the complex plane; thus, it is stoic, and
often more difficult than the corresponding non-elementary approach.
\smallskip
I am now about to relate a deeper message of Riemann. His short but great article was
in fact a r\'esum\'e of his extensive work on $\zeta(s)$. Thus, perhaps because
of his plan of publishing the main article later, he
veiled all of his massive computation and condensed an enormous message into
this famous sentence in the r\'esum\'e:
$$
\eqalignno{
&\hbox{\it Es ist sehr warhscheinlich, dass alle Wurzeln reell
sind.}\cr
&\hbox{\rm It is very probable that all complex zeros are on the 
straight line $\Re s={1\over2}$.}&(\beta) 
}
$$
This $(\beta)$ is the famous Riemann Hypothesis or RH for short; it is 
a mathematical translation of his original, and `reell sind' corresponds to `are
on the  straight line $\Re s={1\over2}\,$'. Again you do not need to understand
$(\beta)$ very correctly, but you are required only to know that Riemann 
surmised that those complex numbers on the right side of $(\alpha)$ should be all
on a particular straight line on the complex plane where $s$ moves around. Indeed
this gives rise to a spectacular view. Thus, on one side we have the real line where
prime numbers
$2,3,5,\ldots, 1481,\ldots$  form a queue as Babylonians and Greeks saw, and on the
other side is a straight line on which are mysterious points called complex zeros of
$\zeta(s)$. Any of the two triggers resonances in the other, via 
the function $\zeta(s)$.
\par
There are analogues of $\zeta(s)$ with which the statements corresponding to
$(\beta)$ hold, but so far no relation has been detected between $(\beta)$ and
those analogues. Huge machine computations have been conducted,
and all the results so far support $(\beta)$; hence most specialists
believe that $(\beta)$ be correct. However, as you may comfortably agree by now,
those numerical data can never be a proof of RH  even if they are awfully
impressive. We are again facing the infinitude the invincible; and
the resolution of Riemann's riddle $(\beta)$ can take place only in our fathomless
minds, if it is ever possible. For further details you are referred to my 
article `Riemann Hypothesis', which is to appear in Sugaku Expositions AMS. I
concluded it with this:  When the deep riddle has ceased to be a great
enigma, how will people gaze the grand event?  Some may monologue, ``From the eternal
distance,  descending gracefully and resplendently, a constellation toward me!'' and
some may sigh out, ``How stupid we were!'' The former is from the last verse of the
poet N, and the latter is a word of the mathematician B.
\smallskip
Returning to $\pi(x)$, the hypothesis $(\beta)$ is known to yield a fine
bound for the error, that is,  $|\pi(x)-\r{li}(x)|<C\sqrt{x}\log x$
with a certain number $C$. Namely, the relative error decreases rather
rapidly as $x$ tends to infinity; a professional expression for the same is 
$$
\pi(x)=\r{li}(x)+O\left(\sqrt{x}\log x\right).\eqno(\gamma)
$$
It is known also that this assertion can be reversed, that is, $(\gamma)$
implies $(\beta)$. Hence, RH is important from the practical point of view as well. 
The strength of RH such as $(\gamma)$ indicates has been advocated by many specialists
from various angles, and by most
popular expositions as well. Obviously I do not need to enhance further such a
prevalent opinion. Thus, I shall rather take an ironical stance, and indicate instead
that RH might not be so strong as is awed generally.  An appropriate example is
offered by another conjecture on primes by Gauss. This concerns the existence of
primes between consecutive squares. Thus, Gauss conjectured that between $n^2$ and
$(n+1)^2$ ($n=1,2,3,\ldots$) should always exist a prime. After a few steps of
experiment, you will be inclined to accept this conjecture. It is,
however, the case that such a seemingly harmless conjecture as this is beyond RH, as
far as the past investigations indicate. This is just a typical instance among
awfully many conjectures on primes, all of which appear unconquerable with RH or with
stronger hypothesis on complex zeros of $\zeta(s)$ either. It appears to me that
$(\beta)$ is a mere tip of an iceberg drifting in a deep bay of a terra
incognita,  which I am sometimes afraid is beyond the reach of our intelligence.
\medskip
Here is a digression about machine computing. I expressed my
cautious opinion already; and this is to endorse it.  Thus, there is or rather
`was' yet another conjecture by Gauss on primes. It asserts that 
$\pi(x)<\r{li}(x)$ for all $x$. However, this conjecture was negated by
J.E. Littlewood by using $(\alpha)$. In spite of this rigorous fact, Gauss'
statement has been found valid for any $x$ with which $\pi(x)$ 
is evaluated. This might sound strange, but it should become natural 
if you learn that the current estimation of the first $x$ where the 
opposite inequality  $\pi(x)>\r{li}(x)$ holds is so huge that confirmation is beyond
even the power of the computers to come in the great distant future. This
number is called the Skews Number. It should be stressed that
concerning this false conjecture of Gauss we do not anymore face the riddle of the
infinitude but only a finite range of real numbers to check; yet we have enormous
difficulty to fix the location of the first change of the sign of
$\pi(x)-\r{li}(x)$.
\par
This reminds me of a myth, perhaps of the Indian origin and seemingly known
in ancient Japan, in which a
gigantic number is recounted. My version is an uncertain sketch of what
I  read when I was a little boy. 
Thus, at the center of the universe is a rock of the
size of a great island. Once every year a celestial damsel dances
over the rock. The ethereal fringe of her diaphanous robe touches once slightly
and erodes the rock. Losing a few atoms every year, the rock shrinks.
The story defines one `kalpa' as to be the number of years until the 
disappearance of the rock. Conventional computers should need more than one
kalpa years to confirm a Skews number. This is a good lesson for us; yes, we
have to use first our brain before relying on machines. In the same
token, I wonder how much we are secured by those large scale 
computations on complex zeros of $\zeta(s)$.
\smallskip
Finally, I shall indicate that RH is hidden even in the
elementary arithmetic that any 10 years old can enjoy.
Thus, let us compute the sum of $1\over2$ and $1\over3$ as to be ${{1+1}\over{2+3}}
={2\over5}$. This is of course absurd. Nevertheless,
let us employ this new way of computation, and see how we can view the world of
fractions or rational numbers. Then, at first we have
$\left\{{0\over1},\,{1\over1}\right\}$, from which we get
$\left\{{0\over1},\, {1\over2},\,{1\over1}\right\}$, 
by this rule of computation. Similarly, we have
$\left\{{0\over1},\,{1\over3},\,{1\over2},\, {2\over3},\,{1\over1}\right\}$. 
Next, we insert those with the denominator equal to $4$, and get
$\left\{{0\over1},\,{1\over4},\,{1\over3},\,{1\over2},\,
{2\over3},\,{3\over4},\,{1\over1}\right\}$. We repeat the same construction, and
find that the steps $5$, $6$, and $7$ are, respectively,
$$
\left\{\txt{{0\over1},\,{1\over5},\,{1\over4},\,{1\over3},\,{2\over5},\,{1\over2},\,
{3\over5},\,{2\over3},\,{3\over4},\,{4\over5},\,{1\over1}}\right\},
$$
$$
\left\{\txt{{0\over1},\,{1\over6},\,{1\over5},\,{1\over4},
\,{1\over3},\,{2\over5},\,{1\over2},\,
{3\over5},\,{2\over3},\,{3\over4},\,{4\over5},\,{5\over6},\,{1\over1}}\right\},
$$
$$
\left\{\txt{{0\over1},\,{1\over7},\,{1\over6},\,{1\over5},\,{1\over4},{2\over7},
\,{1\over3},\,{2\over5},\,{3\over7},\,{1\over2},\,
{3\over5},\,{2\over3},\,{5\over7},\,{3\over4},\,{4\over5},\,{5\over6},\,
{6\over7},\,{1\over1}}\right\}.
$$ 
These sequence of fractions are called
Farey Sequences. If we repeat the procedure indefinitely, they will pack up
between $0$ and $1$. J. Frenel showed that how evenly this packing takes place
is closely related to RH. It may be a surprise that at any class room of 
elementary school is a shadow of RH. You will be amused by watching how the
fractions squeeze in the gaps.
\medskip
Here is a coda with mathematical drums beating to which
you are supposed to abandon yourselves:  I shall show an advanced way to
appreciate those humble sequences. In fact, the notion of the
Farey sequence has played tremendously important r\^oles in the modern number theory.
Thus, as I suggested above,
various  sieve methods are arts of approximations to the rigorism represented by
Eratosthenes' sieve or a stone-heavy method of primality testing of ancient
Alexandria. Also, the divine Riemann zeta-function is a tool to approach to
primes or the Euler product via its $\sum$-definition. The use 
of the latter cannot be perfect because it does not make any real tangent to primes;
that is, it is likewise destined to be an approximative argument. Yes, $\zeta(s)$ is
a  vehicle for us to ride on to come around primes but cannot be an
ultimate gem-detector, unlike non-specialist believe. Hence, we have two main
fronts, united to each other, in the study of the distributions of primes and
integers in general,  in both of which
approximation or compromise, whichever, is inevitable. With this, looking over the
developments in the past, I notice that the whole of Farey sequences is definitely a
central medium for the various applications of the harmonic analysis to the problems
on primes. Namely, in those approximations the decomposition into `wavy' elements is
the philosophy, as in modern physics, and there Farey sequences emerge. This
phenomenon has been observed in a number of modes. Here, 
I shall show briefly the instance that I am most fond of, for 
I encountered it in my own research.
\par
Going to the complex plane with Farey sequences at its fringe, we have the
splendid world of the non-Euclidean geometry and automorphic functions,  which were
both discovered in the 19th century. Thus it should be quite natural to
expect a r\^ole of these concepts to play in the above context. Specifically,
Aufheben via the theory of automorphic functions might be applied to
those approximative arguments involving Farey sequences that were employed
in the theory of $\zeta(s)$. This point of view has been underlying in my research for
more than two decades; and the grand plateau I have recently been to provides a
vista.  To begin the expedition to reach there,  we have to cut a knot
or make it loose at least. Thus, let us pick up pairs of fractions from the last
Farey sequence given above, which are adjacent to each other, and compose the matrices
$$
\cdots,\txt{\left({1\atop3}{2\atop7}\right),\;
\left({2\atop5}{1\atop3}\right),\;
\left({3\atop7}{2\atop5}\right),\;
\left({1\atop2}{3\atop7}\right),\;\left({3\atop5}{1\atop2}\right),\;
\left({2\atop3}{3\atop5}\right),\;\left({5\atop7}{2\atop3}\right),\;
\left({3\atop4}{5\atop7}\right)},\cdots,
$$
with an obvious construction.
These are elements of the unimodular group $\varGamma=\r{PSL}_2(\B{Z})$. 
That is, the world of fractions is covered with the world of unimodular
matrices. On the other hand, $\varGamma$ is a lattice of the Lie group
$G=\r{PSL}_2(\B{R})$. Returning to Farey sequences, they are to fill the 
interval $0$ and $1$ but not completely; thus, we have to handle the matter between
fractions, and here is the r\^ole of harmonic or Fourier analysis.  However,
`between fractions' is too vague, and the argument is destined to become incomplete.
To free ourselves from this difficulty,
we go instead to the fact that elements of Farey sequences, if paired as above,
can be connected with $\varGamma$. Then our view becomes clearer, if not perfect.
Thus, the matter on $G$ that is between elements of $\varGamma$, i.e., inside the
tiles of the tessellation induced by the lattice $\varGamma$ can be well defined  in
terms of $\varGamma\backslash G$. Instead of handling the matter between
fractions, which is indeed near an absurdity, we deal with pairs of fractions,
and gain a portion of relief. Now, into our view comes an infinitely
layered cosmos composed of insurgent $\varGamma$-automorphic  waves
on $G$, that is, $L^2(\varGamma\backslash G)$. 
Recently R.W. Bruggeman and I observed on the plateau that those waves gather
occasionally and tune with a spirit of $\zeta(s)$ (`A new approach to the 
spectral theory of the fourth moment of the Riemann zeta-function'. 
J. reine angew.\ Math., 579 (2005), 75--114). Well, $L^2(\varGamma\backslash G)$
might be a cosmic radio-telescope to observe $\zeta(s)$ a deep-sky object. However,
its resolution power seems inadequate for my purpose. I   
need a more massive telescope.
\bigskip
\noindent
{\csc Conclusion.} Starting with songs or poems, and myths, I have enjoyed
talking about my impression of my field in mathematics. My original
motivation was to introduce you to mathematics; I hope you have got a
glimpse of one of the oldest disciplines that is still dynamic. If it is
fortunately the case, then I would like to suggest you to open your old text books
that you used at elementary school or higher. If you lost them or thew them into
wells like those pupils in Nippur the ancient city of learners, please go to a
bookshop and find new. I assure you that there you will find a fantastic story with
which you must forget time passing, for there you will see that lives of ancient
people are reconstructed with astounding details, better than in any text books on
history. Modern pupils are trying to solve problems, with their backs bending,
in much the same way as the boys and girls of 4000 or more years
ago did. There a universal language is used; yes, so universal
with which you are able to communicate with people of any country, 
any time great past and distant future.
Mathematics can be carried to any place and any time without any exertion
of power, as poems and songs.
\bigskip
\noindent
I thank you very much for your attention.

\vskip 1cm
\noindent
{\csc Yoichi Motohashi}
\smallskip\noindent
Department of Mathematics
\par\noindent
Nihon University
\par\noindent
Surugadai, Tokyo 101--8308, {\csc Japan}
\smallskip\noindent
ymoto@math.cst.nihon-u.ac.jp
\par\noindent
www.math.cst.nihon-u.jp/$\sim$ymoto/
\bye